\documentclass[preprint]{elsarticle}

\usepackage{abstract}
\usepackage{amssymb}
\usepackage{mathrsfs}
\usepackage{amsfonts}
\usepackage{amsmath}
\usepackage{epsfig}
\usepackage{graphicx}
\usepackage{subfigure}
\usepackage{multirow}
\usepackage{bm}
\usepackage[hang,small,bf]{caption}
\usepackage{algorithm}
\usepackage{algpseudocode}
\usepackage{graphicx}
\newcommand*\dashline{\rotatebox[origin=c]{90}{$\dabar@\dabar@\dabar@$}}
\def\Xint#1{\mathchoice 
   {\XXint\displaystyle\textstyle{#1}}%
   {\XXint\textstyle\scriptstyle{#1}}%
   {\XXint\scriptstyle\scriptscriptstyle{#1}}%
   {\XXint\scriptscriptstyle\scriptscriptstyle{#1}}%
   \!\int}
\def\XXint#1#2#3{{\setbox0=\hbox{$#1{#2#3}{\int}$}
     \vcenter{\hbox{$#2#3$}}\kern-.5\wd0}}

\def\dashint{\Xint-}

\title{An Isogeometric Boundary Element Method for elastostatic analysis: 2D implementation aspects}

\author[cardiff]{R.N.~Simpson}
\ead{simpsonr4@cardiff.ac.uk}

\author[cardiff]{S.P.A.~Bordas\corref{cor1}}

\author[cardiff]{H.~Lian}

\author[durham]{J.~Trevelyan}
\ead{jon.trevelyan@durham.ac.uk}

\address[cardiff]{School of Engineering, Institute of Mechanics and Advanced Materials, Cardiff University,
Queen's Buildings,
The Parade,
Cardiff CF24 3AA}

\address[durham]{School of Engineering and Computing Sciences,
Durham University,
South Road,
Durham DH1 3LE}

\cortext[cor1]{Corresponding author}

\date{\today}

\begin{document}


\maketitle
\begin{abstract}
The concept of isogeometric analysis, whereby the parametric
functions that are used to describe CAD geometry are also used to
approximate the unknown fields in a numerical discretisation, has
progressed rapidly in recent years. This paper advances the field
further by outlining an isogeometric Boundary Element Method
(IGABEM) that only requires a representation of the geometry of the
domain for analysis, fitting neatly with the boundary representation
provided completely by CAD. The method circumvents the requirement
to generate a boundary mesh representing a significant step in
reducing the gap between engineering design and analysis. The
current paper focuses on implementation details of 2D IGABEM for
elastostatic analysis with particular attention paid towards the
differences over conventional boundary element implementations.
Examples of Matlab\textregistered\, code are given whenever possible
to aid understanding of the techniques used.
\end{abstract}

\section{Introduction}


Isogeometric analysis (IGA) is a subject which is receiving a great deal of attention amongst the computational mechanics community since it has the potential to have a profound effect on the current engineering design and analysis process. The concept has the capability of leading to large steps forward in efficiency since effectively, the process of ``meshing'' is either eliminated or greatly suppressed.  To understand why this is such a revolutionary development, it is necessary to take an abstract view of the current engineering design and analysis process and evaluate the critical points where bottlenecks occur which subsequently lead to delays in engineering projects.

Models are created in Computer Aided Design (CAD) software allowing
designers to materialise ideas into computational objects that can
range from simple geometries to highly detailed and complex
engineering prototypes. The wide array of CAD packages available,
along with the ever-increasing geometrical modelling capabilities,
offers designers the ability to create realistic models of complex
components. Once a model is complete in CAD,  it must be transformed
into a form suitable for analysis, with the most time-consuming step
taken up in creating a suitable ``analysis-ready'' model which forms
a discretisation of the domain (or boundary).  This step in many
cases requires human intervention by a specialist to ensure that the
discretisation is of sufficient quality to give accurate results in
future simulations. But what is most important to note, is that the
relative portion of time taken to create an analysis-ready design is
approximately 80\% of the total design and analysis process, thereby
dominating the entire design process.

Once a suitable discretisation has been made, analysis can be
carried out using a suitable numerical method such as the Finite
Element Method (FEM), Finite Difference Method (FDM) or the Boundary
Element Method (BEM). Analysis itself represents a relatively small
portion of the design process but importantly, we find that often
the original design is affected by the results obtained from
analysis. In this manner, design and analysis are tightly connected
through an iterative procedure  - a concept which is mirrored
throughout the engineering community.

Recently, an answer to the problems created by the mismatch between
design and analysis was proposed through the concept of Isogeometric
Analysis.  The concept was initially proposed by Hughes et al.
\cite{Hughes:2005p7753} and since this seminal work, a book has been
published entirely on the subject \cite{Cottrell:2009lq}. Rather
than using conventional piecewise polynomial shape functions to
discretise both the geometry and unknown fields, IGA proposes to use
the parametric functions used by CAD as an approximation for both
fields, most commonly taking the form of Non-Uniform Rational
B-Splines (NURBS). In this way, the isoparametric concept is
maintained but more significantly, \emph{the geometry of the problem
is preserved exactly}. In addition, since many of the algorithms
implemented in CAD packages can also be used for numerical analysis,
redundant computations are eliminated allowing analysis to be
carried out with greatly reduced pre-processing.

A great deal of research has been focussed on IGA in recent years
with implementations in areas such as patient-specific modelling
\cite{Bazilevs:2009p8879}, XFEM \cite{Benson:2010p8854}, shells
\cite{Benson:2010p8867} and many others
e.g.\cite{Bazilevs:2008p8895},\cite{Auricchio:2010p8395},\cite{Lu:2011p8841},\cite{Cottrell:2006p8917},\cite{Verhoosel:2011p10082},\cite{Taylor:2011xd}.
In the majority of these methods NURBS are used for discretisation,
but the inability of the functions to produce ``watertight''
geometries and allow local refinement have shown major shortcomings.
Perhaps one of the most significant developments to overcome the
deficiencies of NURBS is the introduction of T-splines
\cite{Bazilevs:2010p8870} and later PHT-splines
\cite{NguyenThanh:2011p8668} that produce watertight geometries and
can be locally refined, benefiting both the design and numerical
analysis communities. In particular,  from an analysis standpoint,
the use of such functions is  essential for efficient algorithms to
exploit adaptivity. More recently IGA has been applied within a BEM
context for elastostatic analysis \cite{Simpson:2011yq} where
particular benefits are realised due to the requirement for only a
surface representation of the geometry. The present paper builds on
this work, where emphasis is given towards the implementation
details of a 2D isogeometric BEM. The organisation of the paper is
as follows: first, an outline of B-splines and NURBS which form the
underlying technology of isogeometric analysis is given, a review of
the conventional BEM is given and finally, the implementation
details of isogeometric BEM are illustrated by building-up an
example problem from an inital CAD model to the final BEM system of
equations.

\section{Geometrical modelling}

The key concept of isogeometric analysis is bringing the fields of design and analysis together into a unified framework through the use of parametric functions that are predominant in CAD. Therefore we concentrate in the current section on describing such functions which most commonly take the form of NURBS. In much of the recent literature on isogeometric analysis \cite{Cottrell:2009lq}\cite{Hughes:2005p7753}\cite{Bazilevs:2010p8870} (and indeed, much literature in the past \cite{Piegl:1995p9119},\cite{Rogers:2001p9055}), extensive details are given on the construction of B-splines and NURBS and therefore we only give a brief description of these functions ensuring that relevant notation is defined to aid the reader in later sections.

Since the present paper is focussed on implementation details, heavy use is made of the algorithms stated in \cite{Piegl:1995p9119} which cover details ranging from evaluating NURBS basis functions to refinement algorithms such as knot insertion and order elevation. The reader is advised to consult this reference since the basics of parametric functions for geometric modelling are outlined clearly.

\subsection{B-splines and NURBS}\label{sec:NURBS}

The first concept which must be grasped when using B-splines or NURBS is that both functions are \emph{parametric} in that the equations which describe the curves (or surfaces) are completely defined by a number of independent parameters. In the current context we will use $\xi$ as the independent variable which describes a B-spline or NURBS and denote it as a coordinate in the \textit{parameter space}. To understand the meaning of the variable $\xi$, an example B-spline is shown in Fig.\ref{fig:figure1}  which illustrates the basic but important concepts of the functions that are used to decribe CAD geometry. We note the following:
\begin{figure}
  \centering
\includegraphics[scale=1]{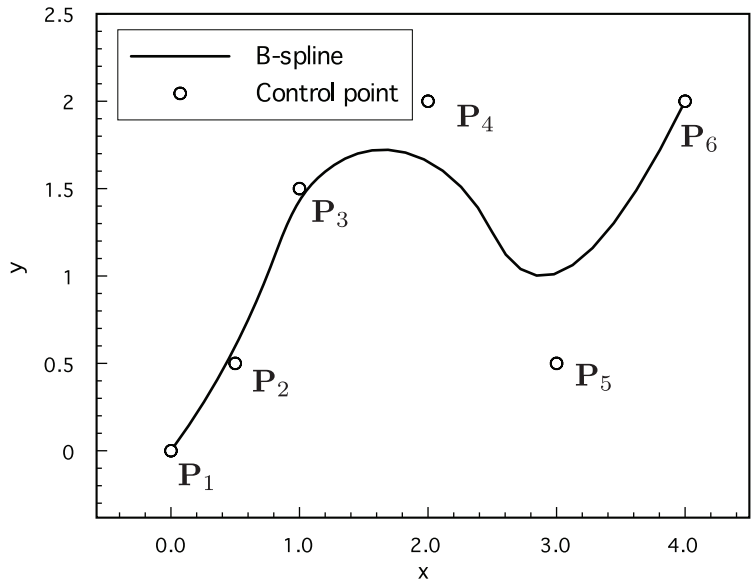}
  \caption{Example B-spline with associated control points.}
  \label{fig:figure1}
\end{figure}

\begin{itemize}
\item The curve requires a set of {\bf control points} to be defined which may or may not lie on the curve
\item The curve requires the definition of a {\bf knot vector}, defined as a non-decreasing sequence of coordinates in the parameter space, which in this case is given by $\Xi = \{0,0,0,1,2,3,4,4,4\}$
\item In this instance the curve is constructed from an \textbf{open knot vector} which results in a spline that is interpolatory at the beginning and end of the curve.
\item If knot values are repeated, then it is found that the order of continuity decreases at that point in the spline.
\end{itemize}


The knot vector is a concept which is often unfamilar to those
working in the field of numerical analysis but should be considered
carefully, since it has a large influence on the resulting spline.
The most important aspect of the knot vector is the relative
difference between the components and for this reason the values can
be scaled if required. That is, the knot vector $\Xi = \{
0,0,0,1,2,3,4,4,4\}$ with $\xi \in [0,4)$ is equivalent to $\Xi = \{
0,0,0,0.1,0.2,0.3,0.4,0.4,0.4\}$ with $\xi \in [0,0.4)$. In later sections we will see that
the process of applying refinement in IGABEM has a considerable
effect on the knot vector. Knot vectors are
ubiqutous in the fields of geometrical modelling and isogeometric
analysis and therefore the reader should become accustomed to their
use.

We now introduce some more formal definitions which are required for
succint notation in later sections.  Denoting the dimensionality of
the problem as $\mathbb{R}^{d}$ ($d=2,3$),
the following three items are required to define fully a B-spline:
\begin{itemize}
\item The curve {\bf degree}, $p$, e.g. linear ($p=1$), quadratic ($p=2$)...
\item A set of of $n$ control points $\mathbf{P}_{a} \in \mathbb{R}^{d}$ , $1\leq a \leq n$
\item A knot vector $\Xi = \{\xi_{1},\, \xi_{2},\,...,\,\xi_{n+p+1}\}$
\end{itemize}

\begin{figure}
  \centering
  \includegraphics[scale=0.3]{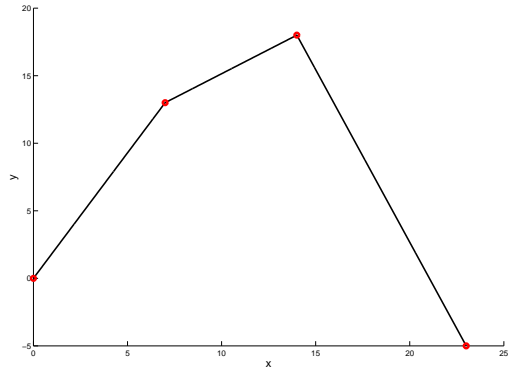}
  \caption{Linear B-spline $p=1$}
  \label{fig:linearBspline}
\end{figure}

\begin{figure}
  \centering
  \includegraphics[scale=0.3]{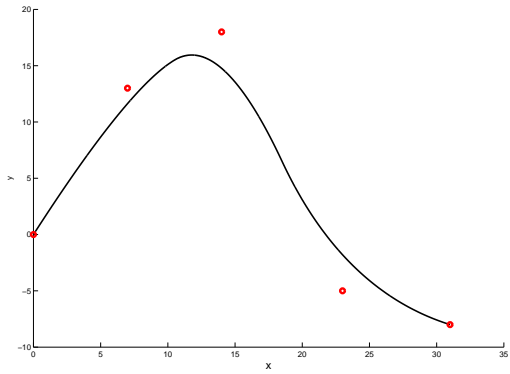}
  \caption{Quadratic B-spline, $p=2$}
  \label{fig:quadBspline}
\end{figure}

\begin{figure}
  \centering
  \includegraphics[scale=0.3]{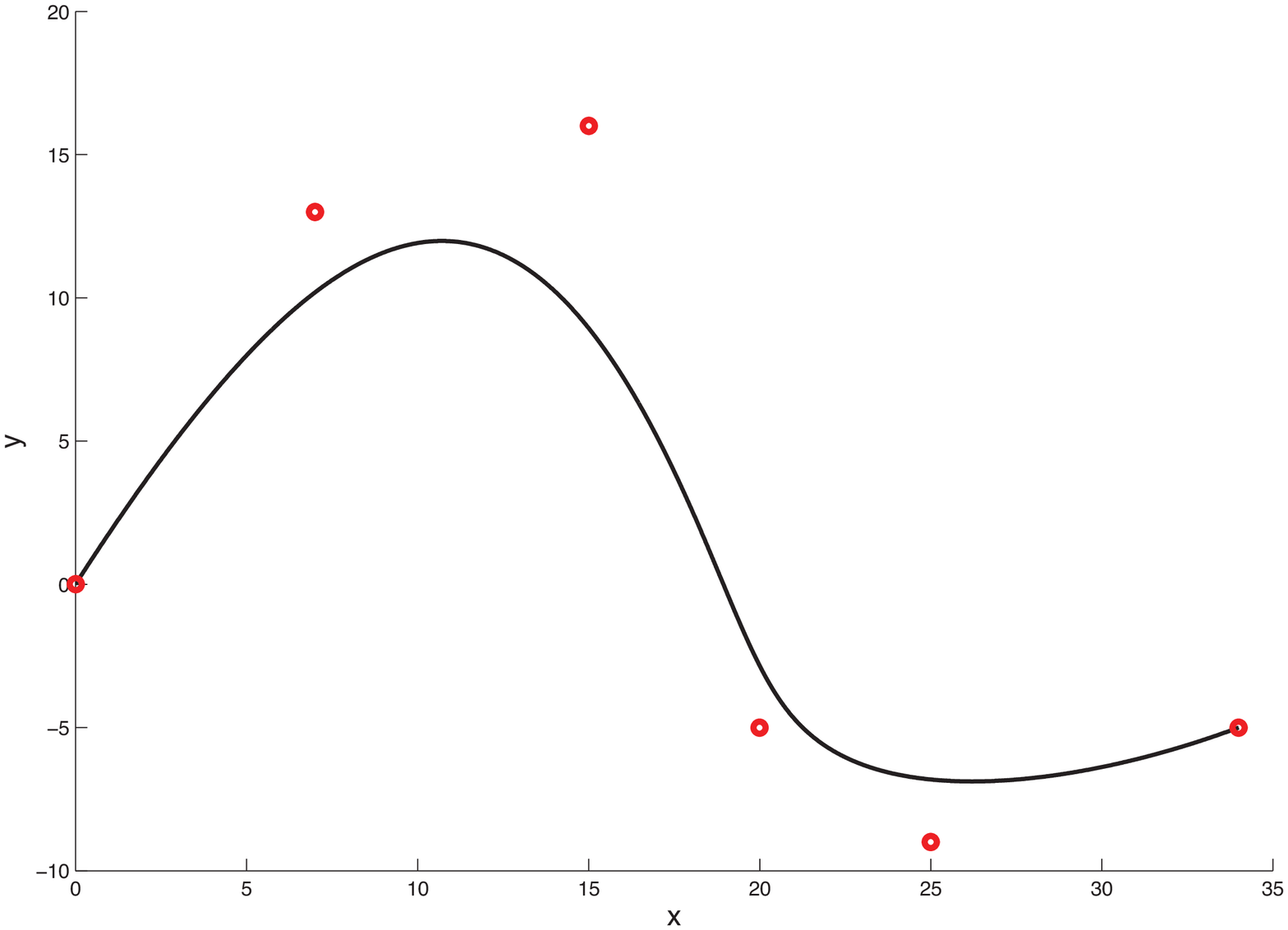}
  \caption{Cubic B-spline, $p=3$}
  \label{fig:cubicBspline}
\end{figure}

Figures~\ref{fig:linearBspline} to \ref{fig:cubicBspline} illustrate
linear, quadratic and cubic B-splines with their associated control
points. The knot vectors associated with each curve are
$\{0,0,1,2,3,3\}$,  $\{0,0,0,1,2,3,3,3\}$ and
$\{0,0,0,0,1,2,3,3,3,3\}$ respectively. Each of these are found to be open knot
vectors which denotes that they contain $p+1$ repeated components at
the beginning and end of the knot vector. The fact that the curve is interpolatory at
the beginning and end is a direct consequence of this. In all future
examples used in this paper, it can be assumed that open knot
vectors are used. Repeated knot components may also occur at points
which are not located at the extremes with a consequence that the continuity of the curve is reduced at that point.


\subsubsection{B-spline and NURBS basis functions}
\label{sec:bsplineBasis}

The previous section served as an overview of B-splines and some of the
terminlogy associated with their construction. But for a B-spline to
be completely defined, some attention must be paid towards their
associated basis functions. The idea of interpolating a discrete
number of points mirrors the technology seen in conventional FEMs
and BEMs but with a distinct difference - the curve is not required
to exhibit the Kronecker delta property at the interpolated points.
The consequences of this are that when interpolating fields such as
displacement, the value obtained at nodal points does not represent
any real displacement but rather a coefficient used for
interpolation. This is similar to meshless methods where techniques
such as Lagrange multipliers \cite{Belytschko:1994oq} and the
penalty method approach \cite{S.N.-Atluri:1998nx} are employed to
impose essential boundary conditions.

Turning our attention towards B-spline basis functions, we introduce the following expression which fully describes a B-spline in terms of its basis functions and control points. This is written as
\begin{equation}\label{eq:bsplineInterpolation}
\mathbf{C}(\xi) = \sum_{a=1}^{n}N_{a,p}(\xi) \mathbf{P}_{a}
\end{equation}
where $\mathbf{C}(\xi)$ is a vector denoting the Cartesian
coordinates of the location described by the parametric coordinate
$\xi$, and $N_{a,p}(\xi)$ denotes the set of B-spline basis
functions of degree $p$ at $\xi$. The basis functions are defined as
\begin{equation}\label{eqn:Ni0}
N_{a,0}(\xi) =
\begin{cases}
1 \quad &\textrm{if} \quad \xi_{a} \leq \xi < \xi_{a+1} \\
0 \quad &\textrm{otherwise}
\end{cases}
\end{equation}
and for $p=1,2,3...$
\begin{align}
N_{a,p}(\xi) = &\frac{\xi - \xi_{a}}{\xi_{a+p} - \xi_{a}} N_{a,p-1}(\xi)\notag\\
&+  \frac{\xi_{a+p+1} - \xi}{\xi_{a+p+1} - \xi_{a+1}} N_{a+1,p-1}(\xi).  \label{eqn:Nip}
\end{align}
\begin{figure}
  \centering
  \includegraphics[scale=0.7]{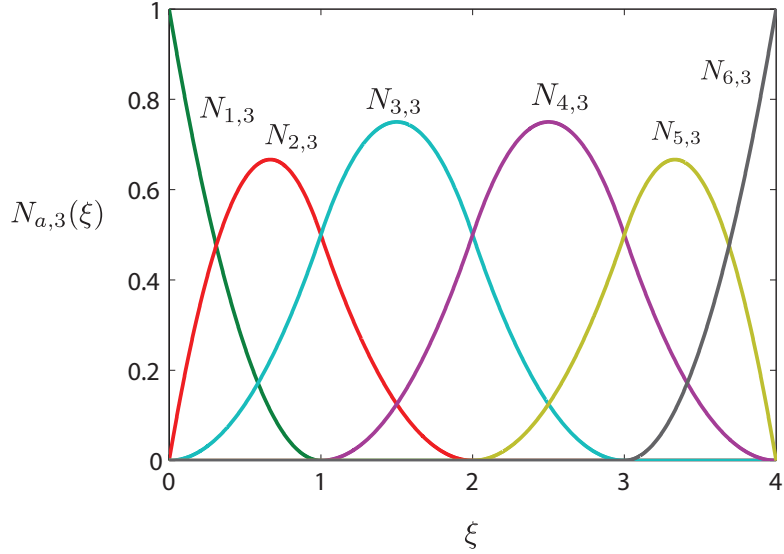}
  \caption{B-spline basis functions for curve shown in Fig.~\ref{fig:cubicBspline} and knot vector $\Xi = \{0,0,0,0,1,2,3,3,3,3\}$}
  \label{fig:bsplinebasis}
\end{figure}

Fig.~\ref{fig:bsplinebasis} illustrates the basis functions
corresponding to the B-spline shown in Fig.~\ref{fig:figure1}; here
the interpolatory nature of the first and last basis functions is
evident. What should be noted is that Eqns.~\eqref{eqn:Ni0} and
\eqref{eqn:NURBSbasis} are recursive in nature and, in their current
form, considerably more expensive than conventional polynomial basis
function expressions. However, there exist several efficient
computational algorithms for their evaluation such as the
Cox-de-Boor algorithm\cite{Piegl:1995p9119},\cite{Rogers:2001p9055}
and, more recently, the extraction operator
\cite{Michael-A.-Scott1:2011ys}. B-spline derivatives are also
required for numerical analysis, but the algorithms for determining
their values are standard in CAD literature, with details given in
\ref{sec:b-spline-derivatives}.

In CAD surface modelling packages, NURBS represent the dominant tool
used to describe curves and surfaces where, in fact, they are found
to be a superset of B-splines and only differ from their
counterparts by the use of an additional coordinate often referred
to as a `weighting'.  In some interpretations NURBS are seen as a
`projection' of B-splines from a higher dimensional space (see
\cite{Cottrell:2009lq}), and it can be shown that some attractive
properties emerge. In particular, NURBS are able to reproduce
circular arcs, spheres and conic sections \emph{exactly} (cf.
B-splines which only approximate such shapes) and this is achieved
through the appropriate choice of weightings.  Each control point
$\mathbf{P}_{a}$ is associated with a weighting $w_{a}$ leading to a
set of NURBS basis functions denoted by $R_{a,p}(\xi)$. The curve is
then interpolated as
\begin{equation}
\mathbf{C}(\xi) = \sum_{a=1}^{n}R_{a,p}(\xi) \mathbf{P}_{a}
\end{equation}
with the NURBS basis functions given by
\begin{equation}\label{eqn:NURBSbasis}
R_{a,p}(\xi) = \frac{N_{a,p}(\xi) w_{a}
}{\sum_{\hat{a}=1}^{n}N_{\hat{a},p}(\xi)w_{\hat{a}}}
\end{equation}
where $N_{a,p}$ can be found from Eqns~\eqref{eqn:Ni0} and
\eqref{eqn:Nip}. In the case that the weights are all set to unity
(i.e. $w_{a}=1 \,\, \forall a$), the basis functions given by
\eqref{eqn:NURBSbasis} reduce to the B-spline basis functions of
\eqref{eqn:Nip}.  Expressions for the derivaties of NURBS basis
functions can also be found, and are detailed in
\ref{sec:nurbs-derivatives}. Since B-splines are a subset of NURBS,
for the sake of generality, NURBS will be considered for all future
examples.

Isogeometric analysis relies on the use of the basis functions outlined in this section and is conceptually very simple - we use the basis functions used to describe the geometry of the problem to approximate the unknown fields in the governing PDEs. That is, in the case of elastostatic analysis using BEM, the displacement and traction components are approximated using NURBS basis functions. The benefit of this approach is clear, since the task of producing a boundary discretisation (mesh) is completely provided by CAD and, once boundary conditions and material properties have been defined, analysis can be carried out immediately.

\section{Conventional BEM}\label{sec:convBEM}

To illustrate the differences between conventional and isogeometric
BEM implementations some details of standard BEM technology are
presented here. The purpose of this section is not to give a
complete derivation of the method but rather an overview to aid
understanding in later sections.  The interested reader is advised
to consult standard BEM references for more details
\cite{Brebbia:1992kb},\cite{Aliabadi:2002qf},\cite{Cruse:1969p9511},\cite{Watson:1979jl}.
In this section we describe the direct collocation form of the BEM
using piecewise polynomial shape functions; the indirect form of the
BEM and the Galerkin BEM are not described, though there appears to
be no reason why the Galerkin BEM cannot be developed in an
isogeometric framework.

To begin, we define the domain of the problem $\Omega$ with boundary
$\Gamma=\partial \Omega$. We also define two points $\mathbf{x'}$
and $\mathbf{x}$, commonly referred to as the source point and field
point respectively, these points being separated by a distance $r$
given by the Euclidean norm
\begin{figure}
  \centering
  \includegraphics[scale=1.0]{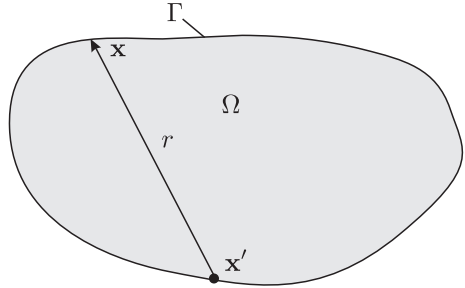}
  \caption{Definition of problem domain with source and field points.}
  \label{fig:bemdomain}
\end{figure}
\begin{equation}
r := ||\mathbf{x'} - \mathbf{x}||
\end{equation}
(see Fig.~\ref{fig:bemdomain}). The point $\mathbf{x'}$ is often
referred to as the collocation point, since in the conventional
collocation BEM implementation, the system of equations is
constructed by taking the collocation point to lie at each nodal
point in turn. The field point $\mathbf{x}$ represents any sampling
point, considered in a numerical integration scheme, on the portion
of boundary over which integration is performed. Making the
assumption of linear elasticity and in the absence of body forces,
we can write the displacement boundary integral equation (DBIE)
which relates displacements and tractions around the boundary
$\Gamma$,
\begin{align}
C_{ij}(\mathbf{x}')&u_{j}(\mathbf{x}') + \dashint_{\Gamma}T_{ij}(\mathbf{x}', \mathbf{x})u_{j}(\mathbf{x})\, \mathrm{d}\Gamma(\mathbf{x})\notag\\
&= \int_{\Gamma}U_{ij}(\mathbf{x}', \mathbf{x})t_{j}(\mathbf{x})\, \mathrm{d}\Gamma(\mathbf{x})\quad i,j=1,2\label{eqn:DBIE}
\end{align}
where $C_{ij}$ is a jump term that arises from the limiting process
of the boundary integral on the left hand side of (\ref{eqn:DBIE})
and is dependent on the geometry at the source point, $u_{j}$ and
$t_{j}$ are the components of displacement and traction around the
boundary and $U_{ij}$ and $T_{ij}$ are displacement and traction
fundamental solutions relating to a source point direction component
$i$ and field point component $j$. These fundamental solutions for
2D linear elasticity may be found in Appendix \ref{sec:2dfundsolns}.

In its current form, Eq.~\eqref{eqn:DBIE} is not amenable for
numerical implementation since $u_{j}$ and $t_{j}$ represent unknown
continuous fields. We therefore proceed in the usual fashion of
discretisation by splitting the boundary of the problem into
elements with local coordinate $\eta \in [-1,1]$ over which the
geometry and fields can be approximated as
\begin{align}
\mathbf{x}_{e}(\eta) &= \sum_{b=1}^{n_{b}} N_{b}(\eta)  \mathbf{x}_{b}\label{eqn:polydiscGeom}\\
\mathbf{u}_{e}(\eta) &= \sum_{b=1}^{n_{b}} N_{b}(\eta) \mathbf{u}_{b}\label{eqn:polydiscU}\\
\mathbf{t}_{e}(\eta) &= \sum_{b=1}^{n_{b}} N_{b}(\eta) \mathbf{t}_{b}\label{eqn:polydiscT}
\end{align}
\begin{figure}
  \centering
  \includegraphics[scale=0.7]{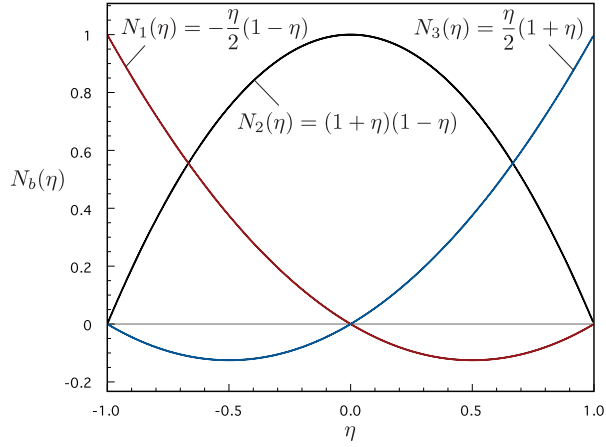}
  \caption{Continuous quadratic basis functions}
  \label{fig:quadbasis}
\end{figure}
where $n_{b}$ is the number of local basis functions (eg. $b=3$ for
a quadratic element), $N_{b}(\eta)$ are the set of polynomial basis
functions (see Fig.~\ref{fig:quadbasis} for the commonly used
quadratic Lagrangian basis functions) and $\mathbf{x}_b$,
$\mathbf{u}_b$, $\mathbf{t}_b$ are vectors of nodal coordinates,
displacements and tractions respectively. The subscript $e$ has been
used in Eqns~\eqref{eqn:polydiscGeom} to \eqref{eqn:polydiscT} to
denote that the vectors apply to a specific element $e$. By
inserting Eqns \eqref{eqn:polydiscU} and \eqref{eqn:polydiscT} into
\eqref{eqn:DBIE}, the discretised DBIE can be written as
\begin{align}\label{eqn:chap3_DiscDBIEConv}
&C_{ij}(\mathbf{x'})u_{j}(\mathbf{x'})\notag\\
&+ \sum_{e=1}^{N_{e}} \sum_{b=1}^{n_{b}} \left[ \int_{-1}^{+1} T_{ij}(\mathbf{x'},\mathbf{x}(\eta))N_{b}(\eta)J^{e}(\eta)\,\mathrm{d}\eta \right]\, u_{j}^{eb}\nonumber\\
\quad &=\sum_{e=1}^{N_{e}} \sum_{b=1}^{n_{b}} \left[ \int_{-1}^{+1}U_{ij}(\mathbf{x'},\mathbf{x}(\eta))N_{b}(\eta)J^{e}(\eta)\,\mathrm{d}\eta \right] \,t_{j}^{eb}.
\end{align}
where $J^e(\eta)$ represents the Jacobian of transformation for
element $e$ that maps $\eta \to \Gamma$ and $1 \leq e \leq N_{e}$ is
the set of element numbers.

As mentioned previously, the system of equations is formed by
considering the collocation point $\mathbf{x'}$ to lie at each nodal
point in turn. In this way, a set of matrices are assembled relating
all displacement components and traction components. This set of
equations can be written as
\begin{equation}\label{eqn:HuGt}
\mathbf{H}\mathbf{u} = \mathbf{G}\mathbf{t}
\end{equation}
with the square matrix $\mathbf{H}$ containing all integrals of the
$T_{ij}$ kernel plus the jump terms $C_{ij}$, the rectangular matrix
$\mathbf{G}$ containing all integrals of the $U_{ij}$ kernel, and
the vectors $\mathbf{u}$, $\mathbf{t}$ containing nodal displacement
and traction components respectively. By prescribing a suitable set
of boundary conditions, which may consist of a set of displacements and
 tractions, equation (\ref{eqn:HuGt}) may be rearranged to give
a set of linear equations in the form
\begin{equation}\label{eq:Axb}
\mathbf{A}\mathbf{x} = \mathbf{b},
\end{equation}
where $\mathbf{x}$ represents the vector of unknown degrees of
freedom. The linear system \eqref{eq:Axb} can be solved
for $\mathbf{x}$ using conventional or fast solvers
\cite{Benedetti:2008p10487,Popov:2001p10517} while noting that
$\mathbf{A}$ is a full, non-symmetric matrix.

\section{Isogeometric BEM}\label{sec:IGABEM}

Our attention now focuses on the main idea of the paper: presenting
the isogeometric boundary element method implementation for 2D
elastostatic analysis.  However, before more details are given, a
comment should be made on the key difference of IGABEM over
conventional BEM. Essentially, it can be reduced to the use of NURBS
basis functions in place of the conventional polynomial
counterparts.  There are certain consequences for  implementation
when NURBS basis functions are used (such as dealing with nodal
points which no longer lie on the boundary), but this key concept
should be kept in mind throughout the following sections.

\begin{figure}
  \centering
  \includegraphics[scale=0.6]{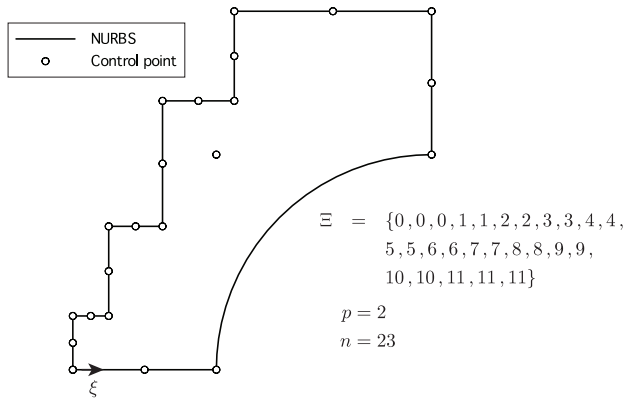}
  \caption{Control points and NURBS curve definition for reactor problem. Appendix~\ref{sec:reactorproblemdata} details the control point coordinates and weights.}
  \label{fig:controlpt_el_defn}
\end{figure}

\begin{figure}
  \centering
  \includegraphics[scale=0.6]{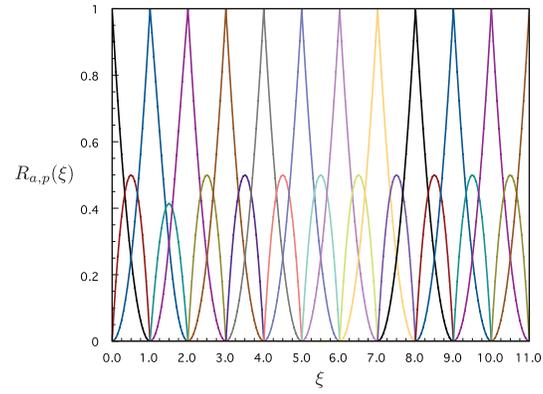}
  \caption{NURBS basis functions for reactor problem.}
  \label{fig:orignalBasisFns}
\end{figure}

We begin by considering a simple 2D geometry of a nuclear reactor
vessel; symmetry has been exploited to simplify the problem. The
exact geometry can be described using the NURBS as illustrated in
Fig.~\ref{fig:controlpt_el_defn}. All information necessary to
define the NURBS curve is provided by a CAD model. In this example
we have chosen, for the entire boundary, to use quadratic basis
functions ($p=2$) which are illustrated in
Fig.~\ref{fig:orignalBasisFns} . 

\subsection{Element definition}
\label{sec:elementdefn}

\begin{figure}
  \centering
  \includegraphics[scale=0.6]{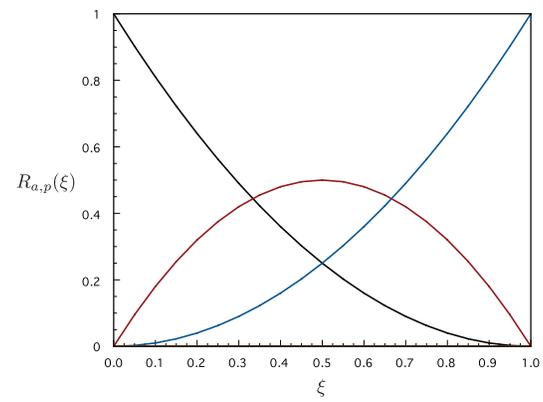}
  \caption{NURBS basis functions for first element of reactor problem. }
  \label{fig:reactorBasis_firstElement}
\end{figure}

As shown in Sec.~\ref{sec:convBEM}, the boundary integrals in the
DBIE are evaluated over the entire domain by summing all elemental
contributions. But for IGABEM, it is not immediately obvious what
our definition of elements should be. However, we should bear in
mind that some notion of `elements' is used in IGABEM simply as a
construct for numerical integration, and this is the definition used
in the remainder of this paper. The element domain is required to
cover only the portion of the boundary where the relevant basis
functions are non-zero. For convenience this leads us to a
definition of element boundaries as the unique values of the knot
vector, which can be implemented simply as
\begin{verbatim}
uniqueKnots = unique(knotVec);
elRanges = [uniqueKnots(1:end-1)' uniqueKnots(2:end)'];
\end{verbatim}
For example, the first element for the reactor problem is defined by
$\xi \in [0,1]$ with the set of non-zero basis functions over this
element illustrated in Fig.~\ref{fig:reactorBasis_firstElement}.

For the purposes of numerical integration using Gauss-Legendre
quadrature, local coordinates in the range $[-1,1]$ over a
particular element are most commonly used, so a transformation that
maps from the parameter space $\xi \in [\xi_{1},\xi_{a+p+1}]$ to a
parent coordinate space defined over an element as $\hat{\xi} \in
[-1,1]$ is required. To achieve this, a Jacobian of transformation
is used which comprises two terms combined in the chain rule sense:
a mapping from the physical coordinate space to parameter space
($d\Gamma/d\xi$) and a mapping from parameter space to the local
parent coordinate space ($d\xi/d\hat{\xi}$). This results in the
following Jacobian of transformation:

\begin{equation}\label{eqn:NURBSjacob}
J(\hat{\xi}) = \frac{d\Gamma}{d\hat{\xi}} = \frac{d\Gamma}{d\xi} \frac{d\xi}{d\hat{\xi}}
\end{equation}
Explicit expressions for the derivatives on the right hand side of
Eq.~\eqref{eqn:NURBSjacob} are given in Appendix~\ref{sec:jacobian}.

\subsubsection{Element connectivity }
\label{sec:elconnectivity}

Before approximations of the geometry and unknown fields can be given, the non-zero basis functions must be determined for a particular element, thus forming a connectivity function. If this is done, then a set of local basis functions that are related to the global basis functions can be defined as
\begin{equation}
N_{b}^{e}(\hat{\xi})\equiv R_{a,p}(\xi(\hat{\xi}))
\end{equation}
where the local basis function number $b$, element number $e$ and global basis function number are related by
\begin{equation}
a = conn(e,b)
\end{equation}
where $conn()$ is a connectivity function. The connectivity function
for the reactor problem is given in
Appendix~\ref{sec:reactor_connmatrix}.  Using this definition of
local basis functions, it is now possible to state isogeometric
approximations for the geometry, displacement and traction as
follows:
\begin{align}
\mathbf{x}_{e}(\hat{\xi}) &= \sum_{b=1}^{p+1} N_{b}^{e}(\hat{\xi}) \mathbf{x}_{b}\label{eqn:NURBSdiscGeom}\\
\mathbf{u}_{e}(\hat{\xi}) &= \sum_{b=1}^{p+1} N_{b}^{e}(\hat{\xi}) \mathbf{d}_{b}\label{eqn:NURBSdiscU}\\
\mathbf{t}_{e}(\hat{\xi}) &= \sum_{b=1}^{p+1} N_{b}^{e}(\hat{\xi})
\mathbf{q}_{b}\label{eqn:NURBSdiscT}
\end{align}
where $\mathbf{x}_{b}$, $\mathbf{d}_{b}$ and  $\mathbf{q}_{b}$ are
vectors of the geometric coordinates, displacement coefficients and
traction coefficients respectively, associated with the control
point corresponding to the basis function $b$. It should be noted
that we have used the term \textit{coefficient} since the NURBS
basis functions do not necessarily obey the Kronecker-delta property
(in contrast to the polynomial basis functions shown in
Fig.~\ref{fig:quadbasis}). Therefore, the terms $\mathbf{d}_b$ and
$\mathbf{q}_b$ do not necessarily represent real displacements and
tractions.

Is it is now a simple case of substituting
Eqns~\eqref{eqn:NURBSdiscU} and \eqref{eqn:NURBSdiscT} into the DBIE
of \eqref{eqn:DBIE} while using the Jacobian given by
\eqref{eqn:NURBSjacob}.  This results in the following discretised
equation for IGABEM:
\begin{align}
C_{ij}&(\mathbf{x'}) \sum_{l=1}^{p+1}N_{l}^{\bar{e}}(\hat{\xi}')d_{j}^{l\bar{e}}\notag\\
 &+ \sum_{e=1}^{N_{e}} \sum_{l=1}^{p+1} \left[\int_{-1}^{+1}T_{ij}(\mathbf{x'}, \mathbf{x(\hat{\xi})})N_{l}^{e}(\hat{\xi})J(\hat{\xi})\, \mathrm{d}\hat{\xi} \right]\, d_{j}^{le}\notag\\
 &= \sum_{e=1}^{N_{e}} \sum_{l=1}^{p+1} \left[\int_{-1}^{+1}U_{ij}(\mathbf{x'}, \mathbf{x(\hat{\xi})})N_{l}^{e}(\hat{\xi})J(\hat{\xi})\, \mathrm{d}\hat{\xi} \right]\, q_{j}^{le}\label{eqn:isoBEMdbie}.
\end{align}
where $d_{j}^{le}$ and $q_{j}^{le}$ represent components of the
vectors $\mathbf{d}_b$ and $\mathbf{q}_b$ for the element $e$. We
denote  $\bar{e}$ as the element containing the collocation point
$\mathbf{x}'$ and $\hat{\xi}'$ is the local coordinate of the
collocation point in element $\bar{e}$. The reason for specifiying
these terms concerns the term $u_{j}(\mathbf{x'})$ as seen in
Eq.~\eqref{eqn:chap3_DiscDBIEConv}. In the conventional
implementation where collocation occurs at nodal points, the
Kronecker-delta property of the basis functions ensures that at the
collocation point itself the basis functions are interpolatory.  In
contrast to this, the IGABEM formulation cannot guarantee that this
is true and the displacement must be interpolated as
$u_{j}(\mathbf{x'}) = \sum_{l=1}^{p+1}
N_{l}^{\bar{e}}(\hat{\xi}')d_{j}^{l\bar{e}}$.

\subsection{Collocation point definition}
\label{sec:collocptdefn}

A significant change in IGABEM over conventional BEM is in the
location of collocation points, since the normal practice of
collocation at nodal positions is no longer valid. We can easily see
why this is the case by inspecting the position of the control
points in Fig.~\ref{fig:controlpt_el_defn}; it is evident that there
is one control point in this example that does not lie on the
boundary. Indeed, for most curved boundaries the control points will
not lie on the boundary. Since the control points can be interpreted
as nodes in the IGABEM formulation, this presents a problem since it
is essential that collocation takes place on the boundary $\Gamma$.
To overcome this, we choose to use the Greville abscissae definition
\cite{Greville:1964p8978}, \cite{Johnson:2005p7969} to define the
position of collocation points in parameter space. This is defined
as:
\begin{equation}\label{eqn:Greville}
\xi_{a}' = \frac{ \xi_{a+1}+\xi_{a+2} + \cdots + \xi_{a+p}}{p} \quad a=1,2,\cdots,n
\end{equation}
In Matlab\textregistered\,, this is easily evaluated as:
\begin{verbatim}
collocPts = zeros(1,n);
for i=1:n
   collocPts(i) = sum(knotVec(i+1:i+p)) / p;
end
\end{verbatim}
where \verb#n# denotes the number of control points, \verb#p# is the
curve order and \verb#knotVec# is the knot vector.
\begin{figure}
  \centering
  \includegraphics[scale=0.6]{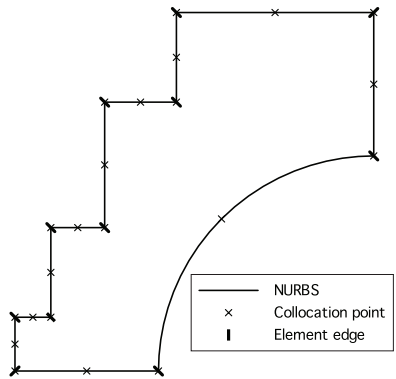}
  \caption{Collocation point and element definitions for reactor problem.}
  \label{fig:collptEldefnReactor}
\end{figure}
If this definition is applied to the geometry of the reactor vessel,
then the collocation points are as shown in
Fig.~\ref{fig:collptEldefnReactor}. The coordinates in physical
space can be found by using \eqref{eq:bsplineInterpolation} with the
NURBS basis functions given by \eqref{eqn:NURBSbasis}. The
definition of the element boundaries is also illustrated in this
figure.


\subsection{Implementation aspects}

The previous two sections outlining the definition of boundary
elements and collocation points within an isogeometric BEM framework
specify the key aspects that must be changed for an existing
collocation BEM implementation. Our attention now focuses on the
implementation for the reactor vessel problem to illustrate these
aspects.

\subsubsection{Refinement}
\label{sec:preprocessing}

The discretisation in Fig.~\ref{fig:controlpt_el_defn} represents
the geometry of the reactor problem exactly, but experienced
analysts will be acutely aware that, although the geometry may be
captured, the basis functions may be insufficient to capture the
gradients in the unknown fields, consequently leading to large
errors in the solution. This leads to one of the most beneficial
properties of B-splines and NURBS for numerical analysis, which is
that the mesh can be refined to arrive at a richer set of basis
functions while \emph{preserving the exact geometry at all stages}.
In the current paper we present two types of refinement: knot
insertion and order elevation, termed h-refinement and p-refinement
in numerical analysis literature.  The algorithms for these
refinement processes are standard in CAD, with details given in
\cite{Piegl:1995p9119},\cite{Rogers:2001p9055} and source code given
in \cite{bobbiesimpsonIsoBEM}.

\begin{figure}
  \centering
  \includegraphics[scale=0.6]{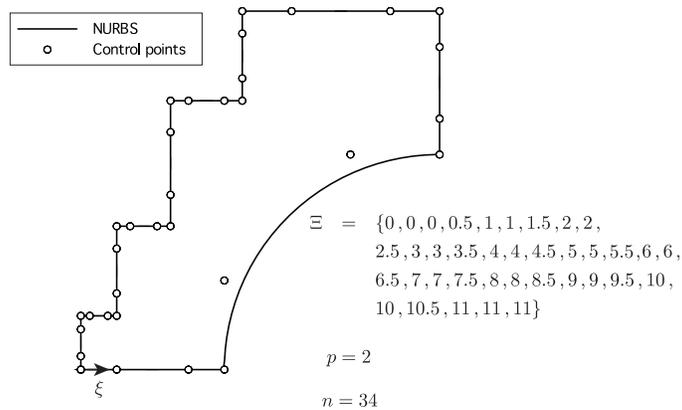}
  \caption{Control points and NURBS curve definition for reactor problem after knot insertion.}
  \label{fig:geomHrefinement}
\end{figure}

\begin{figure}
  \centering
  \includegraphics[scale=0.6]{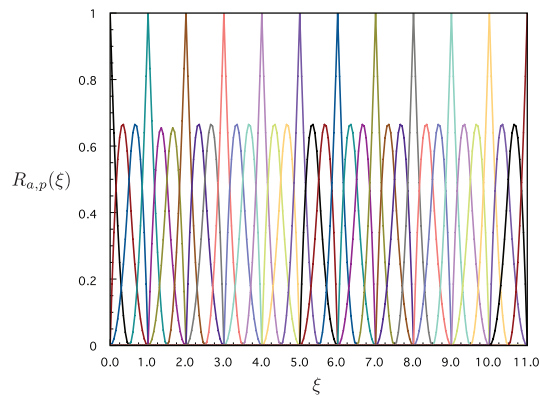}
  \caption{NURBS basis functions for reactor problem after uniform knot insertion (h-refinment)}
  \label{fig:hrefinmentBasis}
\end{figure}

\begin{figure}
  \centering
  \includegraphics[scale=0.6]{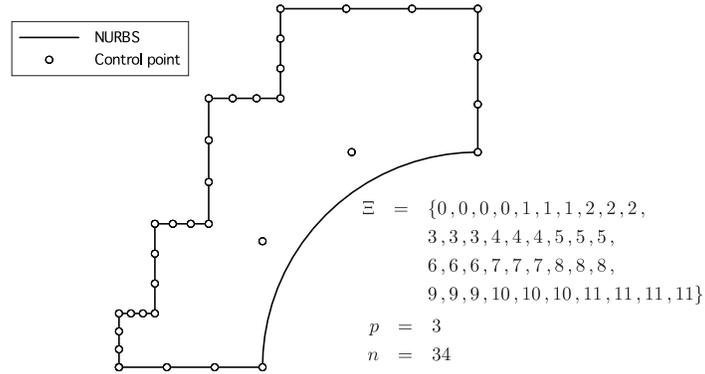}
  \caption{Control points and NURBS curve for reactor problem after order elevation}
  \label{fig:geomPrefinements}
\end{figure}

\begin{figure}
  \centering
\includegraphics[scale=0.6]{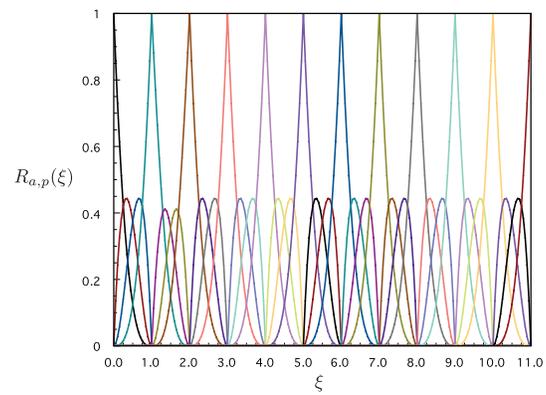}
  \caption{NURBS basis functions for reactor problem after order elevation}
  \label{fig:prefinementBasis}
\end{figure}

We can illustrate both knot insertion and order elevation using the reactor problem shown in Fig.~\ref{fig:controlpt_el_defn}. Figs~\ref{fig:geomHrefinement} and \ref{fig:hrefinmentBasis} illustrate knot insertion where additional knots have been inserted uniformly into the original knot vector. Figs~\ref{fig:geomPrefinements} and \ref{fig:prefinementBasis} illustrate order elevation where the basis functions have been increased from quadratic ($p=2$) to cubic ($p=3$). Both types of refinement introduce changes to the knot vector and set of control points. A third type of refinement which is often referred to as k-refinement can also be considered which is a combination of p- and h- refinement. It arises due to the non-commutative nature of the aforementioned refinements. k-refinement is not considered in the present paper.


\subsubsection{Integration}
\label{sec:integration}

A key feature of any BEM implementation is the evaluation of the
boundary integrals containing the kernels over element domains. It
is well-known that both regular and singular integrands are found
depending on the position of the collocation point relative to the
field element. Essentially, the evaluation of BEM integrals is split
into three different types described as
%
\begin{enumerate}
\item {\bf Regular integration}: the collocation point lies in an element different from the field element.
\item {\bf Nearly singular integration}: the collocation point lies in a element not on but near the field element.
\item {\bf Singular integration}: the collocation point lies in the field element and can be one of two types:
\subitem {\bf Strongly singular} integrals: ($T_{ij}$ kernel,
$\mathcal{O}(1/r)$ in 2D) \subitem {\bf Weakly singular} integrals:
($U_{ij}$ kernel, $\mathcal{O}(\ln(1/r))$ in 2D)
\end{enumerate}
In the present study, we choose to treat the regular and nearly-singular integrals in the same manner, although several methods exists for the efficient treatment of integrals of the latter type \cite{Telles:1987dz},\cite{Lachat:1976hi},\cite{Cruse:1993wf}. The evaluation of singular integrals must be given close consideration since they are found to have a large influence on the accuracy of the resulting solution.

The present work uses the subtraction of singularity method (SST)
\cite{Guiggiani:1987fc} to evaluate strongly singular integrals, in
which the integrand is split into its regular and singular parts.
Further details of the method are shown in \cite{:1998hl} with a
full implementation given in \cite{bobbiesimpsonIsoBEM}. The idea is
to separate the integral into a regular part (which can evaluated
using standard quadrature routines) and a singular part which can be
evaluated analytically (for 2D boundary integrals). If the SST
method is used, however, the jump term $C_{ij}$ must be calculated
explicitly (c.f. the rigid body motion technique
\cite{Cruse:1969p9511} which calculates it implicitly). However, a
simple formula exists \cite{Guiggiani:1987fc} for the 2D case which
is restated here as
\begin{equation}
C=\frac{1}{8\pi(1-\bar{\nu})}\left[
\begin{array}{cc}
\begin{array}{l}4(1-\bar{\nu})(\theta_1-\theta_2)\\\phantom{a}+(\sin2\theta_1-\sin2\theta_2)\end{array}&\cos2\theta_2-\cos2\theta_1\\
\cos2\theta_2-\cos2\theta_1&\begin{array}{l}4(1-\bar{\nu})(\theta_1-\theta_2)\\\phantom{a}-(\sin2\theta_1-\sin2\theta_2)\end{array}
\end{array}
\right]
\end{equation}
where $\bar{\nu} = \nu$ for plane strain and $\bar{\nu} = \nu / ( 1 + \nu)$ for plane stress and the angles $\theta_1$, $\theta_2$ are related to the normals of the surface at the collocation point (see Fig.2 in \cite{Guiggiani:1987fc}).

For weakly singular integrals, a variety of techniques are available including specific logarithmic quadrature points and weights \cite{Stroud:1966fy}, but in the present study we choose to use the Telles transformation \cite{Telles:1987dz} which cancels the singularity leaving a regular integrand. This is achieved through the following transformation:
\begin{equation}
  \label{eq:tellesTransform}
\hat{\xi}=\frac{(\gamma-\gamma')^3 +\gamma'(\gamma'^2+3)}{(1+3\gamma'^2)}
\end{equation}
where
\begin{equation}
  \label{eq:gammaDash}
\gamma'= \sqrt[3]{\hat{\xi}'(\hat{\xi}'^2-1) + |\hat{\xi}'^2-1|} + \sqrt[3]{\hat{\xi}'(\hat{\xi}'^2-1) - |\hat{\xi}'^2-1| + \hat{\xi}'}
\end{equation}
$\xi'$ denotes the location of the singularity in the parent space
($\xi' \in [-1,1]$) and  $\gamma$ represents the new integration
variable. Therefore, a Jacobian which transforms from the parent
space $\hat{\xi} \in [-1,1]$ to $\gamma$ is required. This is given
as:
\begin{equation}
  \label{eq:tellesJacob}
  \mathrm{d}\hat{\xi}=\frac{3(\gamma-\gamma')^2}{1+3\gamma'^2}\mathrm{d}\gamma
\end{equation}
Using Eqn. (\ref{eq:tellesTransform}), (\ref{eq:gammaDash}) and (\ref{eq:tellesJacob}), the transformation of the integration can be expressed as:
\begin{equation}
\int_{-1}^{+1}f(\hat{\xi})\mathrm{d}\hat{\xi}=\int_{-1}^{+1}f\left[\frac{(\gamma-\gamma')^3 +\gamma'(\gamma'^2+3)}{(1+3\gamma'^2)}\right]\frac{3(\gamma-\gamma')^2}{1+3\gamma'^2}\mathrm{d}\gamma
\end{equation}
A derivation of the transformation is given completely in
\cite{Telles:1987dz} with a full implementation for the $U_{ij}$
kernels given in \cite{bobbiesimpsonIsoBEM}.

\subsubsection{IGABEM algorithm}
\label{sec:igabem-algorithm}

Now that the integration routines which form the core of IGABEM
implementation have been outlined, we are in a position to give an
overview of the entire IGABEM algorithm, illustrated in Algorithm
\ref{alg:igabem}.

\begin{algorithm}
\caption{IGABEM algorithm}\label{alg:igabem}
\begin{algorithmic}[1]
\State \textbf{Read CAD input data} \Comment{e.g. Control points, knot vector}
\State \textbf{Read material properties, boundary conditions}
\State \textbf{Perform mesh refinement} \Comment{e.g. knot insertion or order elevation}
\\
\For{$c\gets 1,  n_{c}$}\Comment{Loop over collocation pts}
\State $jumpTerm\gets calcJumpTerm(collocNormals(c))$ \Comment{Calculate $C_{ij}$}
\For{$e\gets 1,  n_{el}$}\Comment{Loop over elements}
\State $elementConn \gets glbConn(e)$\Comment{Element connectivity}
\State $elRange \gets ElmtRanges(e)$\Comment{Range of element}
\\
\If{$c \in elRange$}\Comment{Singular integration}
\State $\mathbf{H}_{sub} \gets SSTIntegration(c,e, C_{ij})$ \Comment{SST integration}
\State $\mathbf{G}_{sub} \gets TellesIntegration(c,e)$ \Comment{Telles transformation}
\Else\Comment{Non-singular integration}
\State $[\mathbf{H}_{sub}, \mathbf{G}_{sub}] \gets GaussLegendreQuad(c,e)$
\EndIf
\\
\State $\mathbf{H}(c,elementConn) \gets \mathbf{H}_{sub}$ \Comment{Add submatrix to global $\mathbf{H}$ matrix}
\State $\mathbf{G}(c,elementConn) \gets \mathbf{G}_{sub}$ \Comment{Add submatrix to global $\mathbf{G}$ matrix}
\EndFor
\EndFor
\\ \State \Comment{Apply boundary conditions}
\State $[\mathbf{A},\mathbf{z}] \gets applyBoundConds(\mathbf{H},\mathbf{G})$
\State $\mathbf{x}\gets solve(\mathbf{A},\mathbf{z})$ \Comment{Solve system of equations}
\end{algorithmic}
\end{algorithm}

In this algorithm the $\mathbf{H}$ and $\mathbf{G}$ matrices have
been calculated explicitly before boundary conditions are applied to
arrive at the final system of equations. For efficiency, commercial
BEM implementations often calculate $\mathbf{A}$ and $\mathbf{z}$
directly, thereby making the construction of the $\mathbf{H}$ and
$\mathbf{G}$ matrices redundant.

\subsection{Example}

\begin{figure}
  \centering
  \includegraphics[scale=0.6]{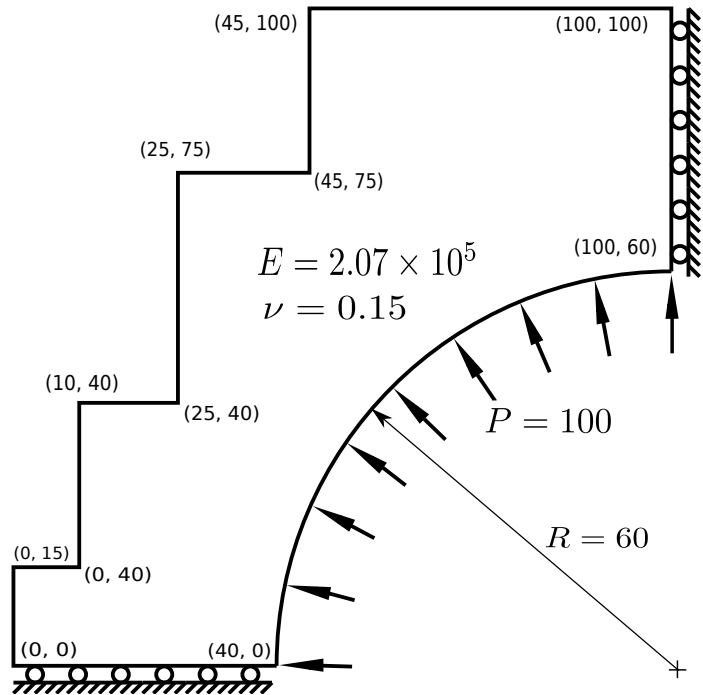}
  \caption{Definition of nuclear reactor geometry, boundary conditions and material properties.}
  \label{fig:reactorDefinition}
\end{figure}

Finally, an example problem using the geometry of the reactor used
throughout the current paper is defined and analysed.  The problem
geometry, boundary conditions and material properties are shown in
Fig.~\ref{fig:reactorDefinition} where symmetrical boundary
conditions are applied at $x=100$ and $y=0$, and a constant pressure
is exerted on the inner boundary. Plane strain is assumed.

\begin{figure}
  \centering
  \includegraphics[scale=0.6]{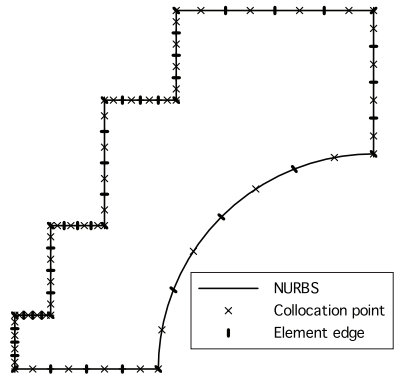}
  \caption{Mesh used for reactor problem IGABEM analysis detailing element edges and collocation point positions.}
  \label{fig:reactorMesh}
\end{figure}

\begin{figure}
  \centering
  \includegraphics[scale=0.6]{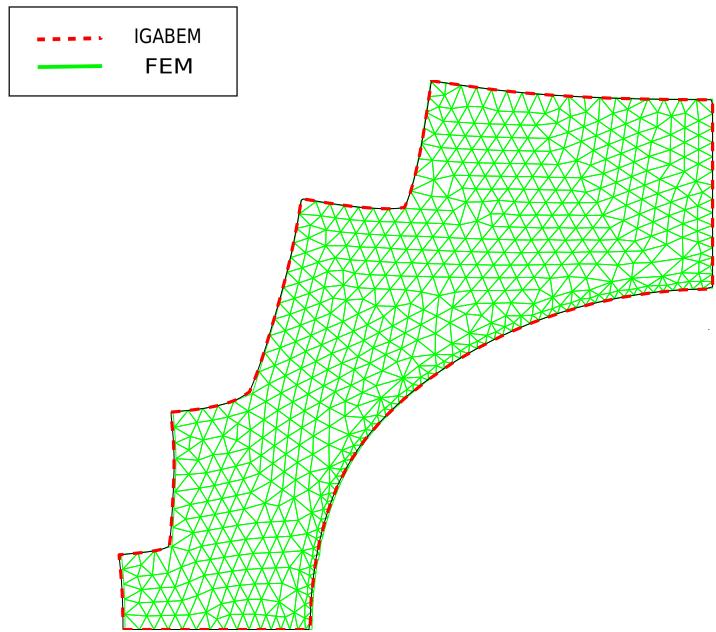}
  \caption{Comparison of IGABEM and FEM results for reactor problem - exaggerated displacement profile.}
  \label{fig:reactorDeformed}
\end{figure}

\begin{figure}
  \centering
  \includegraphics[scale=0.6]{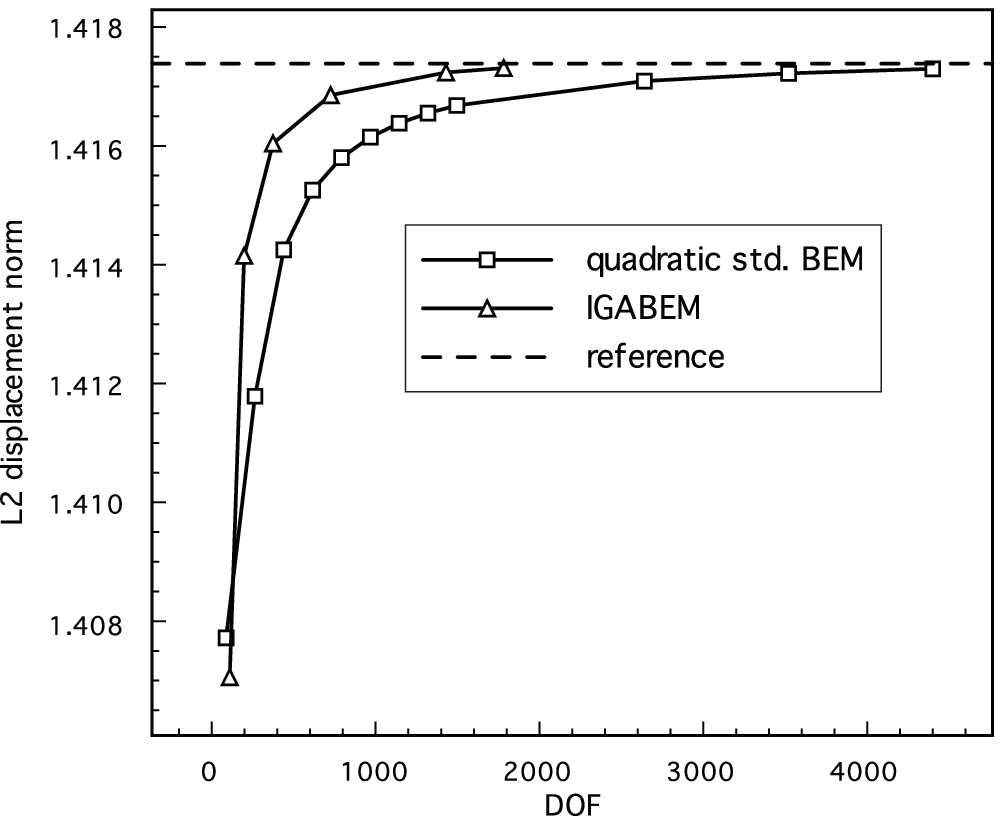}
  \caption{Comparison of $L_2$ displacement norms for reactor problem using IGABEM and standard BEM with quadratic basis functions. }
  \label{fig:L2comparison}
\end{figure}

By applying knot insertion to the original mesh shown in
Fig.~\ref{fig:collptEldefnReactor}, the discretisation shown in
Fig.~\ref{fig:reactorMesh} was used to perform the IGABEM analysis.
The results are shown in Fig.~\ref{fig:reactorDeformed} by plotting
the exaggerated displacement profile with FEM results also plotted
for comparison. As can be seen, excellent agreement is obtained. In
addition, a convergence study was carried out to assess the accuracy
of IGABEM over a standard BEM implementation with quadratic basis
functions. Using h-refinement in both methods, the $L_2$ norm in
displacement was calculated for each mesh as:
\begin{equation}
||\mathbf{u}||_{L_{2}} = \sqrt{\int_{\Gamma} \sum_{i=1}^{d}  (u_{i})^{2}\, \mathrm{d}\Gamma }
\end{equation}
The results obtained are shown in Fig.~\ref{fig:L2comparison} where a
significant improvement in accuracy over the standard BEM
implementation can be seen. The reference solution corresponds to the
converged result of a BEM analysis using quadratic basis functions. 


\section{Conclusions}

The implementation aspects of an isogeometric BEM were outlined,
with attention paid to areas which differ from a conventional BEM
implementation. What is evident is that IGABEM present a
particularly attractive approach for analysis, since the data
provided by CAD can be used directly without the need to create a
mesh. In addition, refinement schemes were outlined which provide
more refinement or richer basis functions in required areas.
Finally, an example was used to illustrate the simplicity and
accuracy of the method where it was shown that good aggrement with FEM was obtained, and in addition, significant improvements in accuracy over a standard BEM implementation with quadratic basis functions were demonstrated.

\appendix

\section{B-splines/NURBS}

\subsection{B-spline derivatives}
\label{sec:b-spline-derivatives}

The first order derivative of the B-spline basis function is
expressed as
\begin{equation}\label{eqn:BsplineDerivatives}
\frac{d}{d\xi}N_{a,p}(\xi)=\frac{p}{\xi_{a+p}-\xi_a}N_{a,p-1}(\xi)-\frac{p}{\xi_{a+p+1}-\xi_{a+1}}N_{a+1,p-1}(\xi).
\end{equation}
where $p$ is the polynomial order, $a$ the basis function index. The higher order derivatives can be obtained by differentiating the two sides of equation (\ref{eqn:BsplineDerivatives}):
\begin{equation}\label{eqn:BsplineHighOrderDerivatives}
\frac{d^k}{d^k\xi}N_{a,p}(\xi)=\frac{p}{\xi_{a+p}-\xi_a}\left(\frac{d^{k-1}}{d^{k-1}\xi}N_{a,p-1}(\xi)\right)-\frac{p}{\xi_{a+p+1}-\xi_{a+1}}\left(\frac{d^{k-1}}{d^{k-1}\xi}N_{a+1,p-1}(\xi)\right).
\end{equation}
From Eqn (\ref{eqn:BsplineDerivatives}) and (\ref{eqn:BsplineHighOrderDerivatives}), we can express high order derivatives with $N_{a,p-k},\cdots,N_{a+k,p-k}$:
\begin{equation}
 \frac{d^k}{d^k\xi}N_{a,p}(\xi)=\frac{p!}{(p-k)!}\sum_{j=0}^{k}\alpha_{k,b}N_{a+b,p-k}(\xi),
\end{equation}
with
\begin{align*}
 \alpha_{0,0}&=1,\\
 \alpha_{k,0}&=\frac{\alpha_{k-1,0}}{\xi_{a+p-k+1}-\xi_a},\\
 \alpha_{k,b}&=\frac{\alpha_{k-1,b}-\alpha_{k-1,b-1}}{\xi_{a+p+b-k+1}-\xi_{a+b}}\qquad  b=1,\dots,k-1,\\
 \alpha_{k,k}&=\frac{-\alpha_{k-1,k-1}}{\xi_{a+p+1}-\xi_{a+k}}.
\end{align*}

\subsection{NURBS derivatives}
From Eqn (\ref{eqn:NURBSbasis}), we can give the first order derivative of NURBS basis function
\label{sec:nurbs-derivatives}
\begin{equation}
\frac{d}{d\xi}R_a^p(\xi)=w_a\frac{W(\xi)N'_{a,p}(\xi)-W'(\xi)N_{a,p}(\xi)}{(W(\xi))^2},
\end{equation}
where $N'_{a,p}(\xi)\equiv\frac{d}{d\xi}N_{a,p}(\xi)$ and
\begin{equation}
 W'(\xi)=\sum_{\hat{a}=1}^{n}N'_{\hat{a},p}(\xi)w_{\hat{a}}.
\end{equation}
 We introduce some notations for convenience
\begin{equation}
 A_{a}^{(k)}(\xi)=w_a\frac{d^k}{d\xi^k}N_{a,p}(\xi), \qquad (\textrm{no sum on}\ i)
\end{equation}
and
\begin{equation}
 W^{(k)}(\xi)=\frac{d^k}{d\xi^k}W(\xi).
\end{equation}
Then higher-order derivatives of these rational functions may be expressed in terms of lower-order derivatives as
\begin{equation}
 \frac{d^k}{d\xi^k}R_a^p(\xi)=\frac{A_a^{(k)}(\xi)-\sum_{b=1}^k\binom{k}{b}W^{(b)}(\xi)\frac{d^{(k-b)}}{d\xi^{(k-b)}}R_a^p(\xi)}{W(\xi)},
\end{equation}
where
\begin{equation}
\binom{k}{b} = \frac{k!}{b!(k-b)!}.
\end{equation}

\section{BEM}

\subsection{2D Fundamental solutions}
\label{sec:2dfundsolns}

Denoting $\mu$ as the shear modulus, $\nu$ as Poisson's ratio, $\delta_{ij}$ as the kronecker-delta function defined by 
\begin{equation}
  \label{eq:krondelta}
\delta_{ij} = 
\begin{cases}
  &0 \quad i\neq j\\
  &1 \quad i=j,
\end{cases}
\end{equation}
and noting that a comma denotes differentiation, the fundamental solutions for 2D linear elasticity are given by:
\begin{equation}
U_{ij}(\mathbf{x}',\mathbf{x})=\frac{1}{8\pi\mu(1-\nu)}\left\{(3-4\nu)\ln\left(\frac{1}{r}\right)\delta_{ij}+r_{,i}r_{,j}\right\}
\end{equation}
\begin{equation}
T_{ij}(\mathbf{x}',\mathbf{x})=\frac{-1}{4\pi(1-\nu)r}\left\{\frac{\partial r}{\partial n}[(1-2\nu)\delta_{ij}+2r_{,i}r_{,j}]-(1-2\nu)(r_{,i}n_j-r_{,j}n_i)\right\}
\end{equation}

\subsection{Boundary element integration parameters}

\subsubsection{Normals}
\label{sec:normals}

The normals can be calculated by:\\
\begin{equation}
n_x=\frac{1}{J(\xi)}\left[\frac{\mathrm{d}y(\xi)}{\mathrm{d}\xi}\right],
\end{equation}
\begin{equation}
n_y=\frac{1}{J(\xi)}\left[\frac{\mathrm{d}x(\xi)}{\mathrm{d}\xi}\right],
\end{equation}
with
\begin{equation}
\frac{\mathrm{d}x(\xi)}{\mathrm{d}\xi}=\sum_{b=1}^{p+1}\frac{\mathrm{d}N_b(\xi)}{\mathrm{d}\xi}x_b,
\end{equation}
\begin{equation}
\frac{\mathrm{d}y(\xi)}{\mathrm{d}\xi}=\sum_{b=1}^{p+1}\frac{\mathrm{d}N_b(\xi)}{\mathrm{d}\xi}y_b.
\end{equation}

\subsubsection{Jacobian of transformation}
\label{sec:jacobian}

\begin{equation}
\frac{\mathrm{d}\xi(\hat{\xi})}{\mathrm{d}\hat{\xi}}=\frac{\xi_f-\xi_s}{2},
\end{equation}
where $\xi_s$ and $\xi_f$ denotes the values of knots at the beginning and end of the element respectively (assuming an outward pointing normal).\\

\begin{equation}
\frac{\mathrm{d}\Gamma}{\mathrm{d}\xi}=\sqrt{\left(\frac{\mathrm{d}x(\xi)}{\mathrm{d}\xi}\right)^2+\left(\frac{\mathrm{d}y(\xi)}{\mathrm{d}\xi}\right)^2}.\\
\end{equation}

\section{Reactor problem data}
\label{sec:reactorproblemdata}

\subsection{Control points and weights}
\label{sec:reactorproblem_cptsweights}

\label{tab:controlPoints}

\begin{tabular}{c|c|c}
\hline
index ($a$) & Control point coordinate ($\mathbf{P}_a$)& Weight $w_a$\\
\hline
1&(0, 0)&1\\
\hline
2&(20, 0)&1\\
\hline
3&(40, 0)&1\\
\hline
4&(40, 60)&$\sqrt{2}/2$\\
\hline
5&(100, 60)&1\\
\hline
6&(100, 80)&1\\
\hline
7&(100, 100)&1\\
\hline
8&(72.5, 100)&1\\
\hline
9&(45, 100)&1\\
\hline
10&(45, 87.5)&1\\
\hline
11&(45, 75)&1\\
\hline
12&(35, 75)&1\\
\hline
13&(25, 75)&1\\
\hline
14&(25, 57.5)&1\\
\hline
15&(25, 40)&1\\
\hline
16&(17.5, 40)&1\\
\hline
17&(10, 40)&1\\
\hline
18&(10, 27.5)&1\\
\hline
19&(10, 15)&1\\
\hline
20&(5, 15)&1\\
\hline
21&(0, 15)&1\\
\hline
22&(0, 7.5)&1\\
\hline
23&(0, 0)&1\\
\hline
\end{tabular}

\subsection{Connectivity matrix}
\label{sec:reactor_connmatrix}

\begin{tabular}{c|c}
\hline
element index ($e$)& global basis index ($a$) for ($b_1$,$b_2$,$b_3$) \\
\hline
1&1, 2, 3\\
\hline
2&3, 4, 5\\
\hline
3&5, 6, 7\\
\hline
4&7, 8, 9\\
\hline
5&9, 10, 11\\
\hline
6&11, 12, 13\\
\hline
7&13, 14, 15\\
\hline
8&15, 16, 17\\
\hline
9&17, 18, 19\\
\hline
10&19, 20, 21\\
\hline
11&21, 22, 1\\
\hline
\end{tabular}

v

\bibliography{localbib}
\bibliographystyle{plain}

\end{document}